\documentclass[12pt]{article}
\title {The Brownian loop soup}
\author {Gregory F. Lawler\footnote {Cornell University; Research supported in part by the 
National Science Foundation}  \and 
Wendelin Werner\footnote {Universit\'e Paris-Sud and IUF}}
\usepackage{graphicx,amssymb}
\newtheorem{corollary}{Corollary}
\newtheorem{theorem}[corollary]{Theorem}
\newtheorem{lemma}[corollary]{Lemma}
\newtheorem{conjecture}{Conjecture}
\newtheorem{proposition}[corollary]{Proposition}
\newcommand{\Prob} {{\bf P}}
\newcommand{\R}{\mathbb{R}}
\newcommand{\C}{\mathbb{C}}

\def \H{\mathbb{H}}
\def\Half {\mathbb{H}}
\def \eps {\epsilon}
\def \P {{\bf P}}

\def \E {{\bf E}}
\def \p {{\partial}}

\def \dist {{\rm dist}}
\def \Im {{\rm Im}}
\def \Re {{\rm Re}}
\def \Disk {{\mathbb D}}

\def \rect {{\cal R}}

\def \Pf {\noindent {\bf Proof.} }
\def\endpf {\hfill $\square$  \medskip}
\def \rad {{\rm rad}}
\def \hCp  {{\rm hcap}} 
\def \hcap {\hCp}

\def \percurves {\tilde {\fincurves}}
\def \unrootedcurves {\percurves_{\rm U}}
\def \unrootedmetric {{d_{p,U}}}
\def \metricmeasure {{\cal M}}
\def \fincurves {{\cal K}}
\def \finmetric {d_\fincurves}

\def \loopmeasure{\mu^{\rm loop}}

\def \area {{\rm area}}
\def \bubble {\mu^{{\rm bub}}}

\def \normed{\#}

\begin{document}
\maketitle

\begin {abstract}
We define a natural conformally invariant measure on unrooted 
Brownian loops in the plane and study some of its properties.
We relate this measure to a measure on loops
rooted at a boundary point of a domain and show how
this relation gives a way
to ``chronologically add Brownian loops'' to simple curves in the plane.
\end {abstract}

\section{Introduction}

The recent study of conformally invariant
scaling limits of two-dimensional
lattice systems has shown that measures on
paths that satisfy conformal invariance
(or conformal covariance) and a certain
restriction property are important.  
In particular, in \cite{LSWrest}, it is shown how to construct
``restriction'' measures by dynamically adding  bubbles to
Schramm-Loewner evolution (SLE) curves. As announced there, this 
construction has an equivalent formulation in terms of a Brownian soup of 
loops. The purpose of this paper is to describe these Brownian measures and 
to prove this equivalence.

This description will be given without reference to SLE and is
interesting on its own, but since this is what initiated our 
own interest, let us now describe the link with SLE.  
In \cite {S}, Oded Schramm introduced the SLE processes. These 
are the only random  non-self-crossing 
curves in a domain that combine conformal invariance and a certain
Markovian-type property.
 The definition of  SLE is based on these two facts and can be 
viewed as a dynamic construction: one constructs the law of the curve on the time-interval $[t, t+dt]$
given $\gamma [0,t]$ and then iterates this procedure.
In \cite {LSWrest}, following ideas of \cite {LW2} and partially motivated by the 
problem of the self-avoiding walks in the plane (see \cite {LSWsaw}), a different approach 
to SLE was described. Basically, one looks at how the law of the random curve (seen globally) is 
distorted by an infinitesimal perturbation of the domain it is defined in.
It turns out that a one-dimensional 
family of random sets is in some sense invariant under such perturbation. 
These are called 
restriction measures in \cite {LSWrest}, where it is shown that all of
these measures are  closely related to 
Brownian excursions.
The law of SLE, except for the special case of the SLE with parameter $
\kappa=8/3$, is  not a restriction measure. 
However, one can measure precisely the ``restriction defect'' (i.e., 
the Radon-Nikodym derivative)
with a term involving the Schwarzian derivatives 
of the corresponding conformal maps. On interpretation
 for SLE$_\kappa$ with $\kappa < 8/3$ 
goes as follows: if one adds a certain Poissonian cloud of Brownian bubbles
to the SLE curve, then 
the resulting set is restriction invariant. This can 
be understood simply when $\kappa=2$. In that case,
the SLE curve is \cite {LSWlesl} the scaling limit of the loop-erased random walk. 
In the scaling limit, the corresponding random walk converges to the Brownian excursion 
(which is restriction invariant). Hence, it is 
not surprising that if one puts the erased Brownian bubbles back onto the SLE$_2$ curve, 
one obtains a restriction 
measure. This Poissonian cloud of Brownian bubbles provides a  simple geometric 
picture of the distortion of the law of SLE  under perturbation of the boundary of a domain.
This ``variational'' approach to SLE is closely related some conformal field theory considerations of e.g. \cite {BPZ, Ca}, as pointed out in \cite {FW,FW2}.
The density of the Poissonian cloud in particular plays the role of the (negative of the)
central charge of the corresponding model in the theoretical physics language.

\medbreak

We will describe various measures on Brownian paths with an emphasis on two measures,
the Brownian loop measure and the Brownian bubble measure. The latter was already 
defined and used in \cite {LSWrest} for the previously described reasons.

The Brownian loop measure  is an infinite measure on unrooted Brownian loops in the plane.
It is defined on the set of periodic continuous functions in the plane, where two functions are considered to be indistinguishable if one is obtained by a simple
translation in time ($t \mapsto t + c$), and we call these equivalence classes
 ``unrooted loops''. The Brownian loop measure is scale invariant, and translation-invariant.
Furthermore, it is conformally invariant in the following sense: If there exists a conformal map $\phi$ from $D$ onto $D'$, then the image under $\phi$ of the Brownian loop measure restricted to those loops that stay in $D$ is exactly the Brownian loop measure restricted to those loops that stay in $D'$. This property is in fact very closely related to the 
restriction property.
 
This measure can also be considered a measure on ``hulls'' (compact sets $K$ such
that $\C \setminus K$ is connected) by ``filling in'' the bounded loops. 
It is possible to argue that the Brownian loop measure 
is the only measure on hulls that is conformally 
invariant in the previous sense. 

The Brownian loop soup of intensity $\lambda >0$
  is a realization of a Poisson point process 
of density $\lambda$ times the Brownian loop measure. In other
 words, a sample of the Brownian loop soup is a 
countable family of Brownian unrooted loops.   
There is no non-intersection condition or other interaction between the 
loops. Each loop will  intersect countably many other loops in the same realization of the loop soup.
Although for some purposes it is sufficient to consider the hull generated by a loop,
we will study the measure on loops
with  time  parametrization in this paper.
This is partially motivated by possible future applications.

\medbreak
A bubble in a domain $D$ will be a continuous path 
$\gamma [0,T]$ such that $\gamma (0,T) \subset D$ and 
$\gamma (0) = \gamma (T) \in \partial D$.
We say that the bubble is rooted at $x$ if $\gamma (0) = x$.
The Brownian bubble measure was introduced in 
\cite {LSWrest} in order to construct the restriction measures via SLE.  
The Brownian bubble in $D$, rooted  in $x \in \partial D$
is a $\sigma$-finite measure on Brownian loops that start
and end at $x$, and otherwise stay in $D$.
The description is simplest if the 
considered domain is the upper half-plane $\H$, and the root is the origin.
We will see that it can be considered
as a conditioned version of the Brownian loop 
measure. The relation between these two measures will lead to 
an equivalence that we now describe.

Loosely speaking, the relation is as follows.  Imagine that
a realization of the loop soup in $\H$ has been chosen, but we cannot see
a loop until we visit a point on that loop.   
Suppose that we travel along a simple 
curve $\eta$ with $\eta (0) = 0$ and $\eta (0, \infty) \subset \H$.
Each time $t$ at which one encounters a loop in the loop soup for the first time, we
can see the whole loop. This prescribes the order in which one finds the loops that
intersect the curve $\eta[0,\infty)$.
 These loops are a priori unrooted; however, we can 
 makes them into rooted loops by starting a loop found at time
$t$ at the point  $\eta(t)$ . If we use this point as a 
root, the loop becomes a bubble in the domain $\H \setminus \eta[0,t]$.
The point is that this loop is ``distributed'' according to the bubble measure. 
    
More precisely, let $\eta$ be as before.
We do not make smoothness assumptions 
on $\eta$; in fact, the cases of most interest to us are
 SLE curves that have Hausdorff dimension greater than one.
Assume that $\eta$ is parametrized by its ``half-plane capacity'' 
(as is customary for Loewner chains in the 
 upper half-plane), i.e.,  that for all $t$,  there exists a (unique) conformal map $\tilde g_t$ from 
 $\H \setminus \eta  [0,t]$ onto $\H$ such that
$\tilde g_t(z) = z + 2t/z  + o(1/z)$ when $z \to \infty$.  We let
$g_t(z) = \tilde g_t(z) - \tilde g_t(\eta_t)$.

Suppose that the countable collection 
of loops $\{\gamma_1,\gamma_2,\ldots\}$ in $\Half$ is a 
realization of the Brownian loop soup in $\H$ with intensity $\lambda
> 0$. This is a random family of equivalence classes
of curves $\gamma_j:[0,t_j] \rightarrow \Half$ with $\gamma_j(0) = \gamma_j(t_j)$, under
the equivalence $\gamma^1 \sim \gamma^2$ if the time-lengths $t^1$ and $t^2$
of $\gamma^1$ and $\gamma^2$ are identical, and if  for some $r$,
$\gamma^1(t) = \gamma^2(t+r)$ for all $t$ (with addition modulo $t^1$). 
For each $j$, let 
\[
r_j = \inf\{s: \eta(s) \in  \gamma_j[0,t_j] \} 
\]
with $r_j = \infty$ if $\eta[0,\infty) \cap \gamma_j[0, t_j] = \emptyset$. 
It is not difficult to see that with probability one for each $j$ with $r_j < \infty$
there is a unique $t \in [0,t_j)$ such that $\gamma_j(t) = \eta(r_j)$.
Then we can choose the representative  $\gamma_j$   so that $\gamma_j(0) = \eta(r_j)$. 
Note that $\gamma_j$ is a bubble in $\H \setminus \gamma [0, r_j]$. We
define $\tilde \gamma_{r_j}$ as the image of 
$\gamma_j$ under the mapping $g_{r_j}^{-1}$ where the
 time-parametrization 
of $\tilde \gamma_{r_j}$ is obtained from that of
$\gamma_j$ using the usual Brownian time-change under conformal maps.
Note that each $\tilde \gamma_{r_j}$ is a bubble rooted at the origin in $\H$. 
Here, the parametrization of $\eta$ by its half-plane capacity is important since we 
index $\tilde \gamma$ by the time $r_j$.

\begin {theorem}
\label {ls=bs}
The process $(\tilde \gamma_r , r \ge 0)$ is a Poisson point process 
with intensity  $\lambda$ times the Brownian bubble measure in $\H$.
\end {theorem}

Of course, this statement depends on the precise definitions of these measures, but it shows that adding the Poisson cloud of bubbles (as in \cite {LSWrest}) to the path $\eta$ is exactly the same as adding to $\eta$ the set of loops in a loop soup that it does intersect. 

One direct application is that adding Brownian bubbles to $\eta$ or to the time-reversal 
of $\eta$ (that is, viewing $\eta$ as a curve from $\infty$ to the origin) is the same 
(from the point of view of the outside hulls).
In the case where $\eta$ is chordal $SLE_2$, it corresponds to the fact that 
loop-erasing a random walk does not depend (in law) on the chosen time-orientation.
This result is more generally closely related to the question of reversibility of the 
SLE's.

One other application is for the ``duality'' conjecture of the SLEs, see \cite {Dub}:
Indeed, from the point of vue of the outside hulls, adding the loops of 
the loop soup to a curve or to its outer boundary is the same. Hence, 
the same is true if one adds (dynamically) bubbles to a curve or to its outer boundary.
This leads to an identity in law between the set obtained by adding the same loop soup to a process closely related to $SLE_\kappa$ or to a process closely related to $SLE_{16/\kappa}$. See \cite {Dub} for more details. 
Theorem~\ref {ls=bs} is also used in  \cite {Wgi}.

Another main point is just the definition of the Brownian loop measure.
Despite its simplicity (and maybe its importance) and its nice properties,
it does not seem (to our knowledge) to
have been considered before. 

The technical aspects of the present paper are not difficult.
Once one has the correct definitions, the proofs are more or less standard exercises on
Brownian motions, excursion theory and Green functions. In order to
 keep  the pace of the paper flowing, we will at times be somewhat informal
(we will not always describe precisely how to take the limit 
of one measure on paths, etc.), leaving the gaps to the interested reader. We will
however  not completely omit these problems (see, e.g.,  the next section).

The paper is organized as follows. In the next section, we mainly introduce some notation. In Section 3, we define some measures on Brownian paths, among which the Brownian bubbles. These are not new, but it is convenient to summarize some of their features in order to 
simplify the relation with the Brownian loop measure. This measure is defined and studied in 
Section 4, the relation with the bubble measure is described in Section 5. The final section is devoted to the question of time-parametrization of the Brownian ``loop-adding'' procedure.

\section{Notations}

We will write $\Disk$ for the unit disk, $\Half = \{x+iy:
y > 0\}$ for the upper half-plane, and $\Disk_+$
for 
$  \Disk \cap \Half = \{z \in \Half:
       |z| < 1 \} . $

 Let $\fincurves$
be the set of all parametrized continuous planar curves 
 $\gamma$ defined on a time-interval $[0,t_\gamma]$.  
 We consider $\fincurves$ as a metric space with the metric
 \begin{equation}  \label{jan20.1}
  \finmetric(\gamma,\gamma^1) =
       \inf_\theta \; [\; \sup_{0 \leq s \leq t_\gamma}
      |s - \theta(s) |  + |\gamma(s) - \gamma^1(\theta(s))|  
            \; ], 
\end{equation}
where the infimum is over all increasing 
homeomorphisms $\theta: [0,t_\gamma] \rightarrow [0,t_{\gamma^1}]$.  
Note that $\fincurves$ under this metric does not identify curves
that are the same
  modulo time-reparametrization.

If $\mu$ is any measure on $\fincurves$, we let
$|\mu| = \mu(\fincurves)$ denote the total mass.
If $0 < |\mu| < \infty$, then we let $\mu^\normed
= \mu/|\mu|$ be $\mu$ normalized to be a probability
measure. 

Let $\metricmeasure$ denote the set of finite
Borel measures on $\fincurves$.  This is a metric
space under the Prohorov metric $d$ (see \cite[Appendix III] {Billingsley}, e.g., for details).
When we say that a sequence of measures converges it will
be with respect to this metric. 
Recall that one standard way to show that two
probability measures $\mu$ and $\nu$ are close with respect to this metric is via 
coupling: 
one finds a probability measure $m$ on $\fincurves \times \fincurves$
 whose first marginal is $\mu$, whose second marginal is $\nu$, and such that
\[   m [\{(\gamma^1,\gamma^2): \finmetric(\gamma^1,\gamma^2)
  > \epsilon\}] \leq \epsilon . \]
To show that a sequence of finite
measures $\mu_n$ converges to a finite measure $\mu$,
it suffices to show that $|\mu_n| \rightarrow |\mu|$ and $\mu_n^\# \rightarrow
\mu^\#.$

If $D$ is a domain, we say that $\gamma$ is 
in $D$ if $\gamma(0,t_\gamma) \subset D$; note that we do not require the endpoints 
of $\gamma$ to be in $D$. Let $\fincurves(D)$ be the set of $\gamma \in \fincurves$ that are in $D$.  
If $z,w \in \C$, let $\fincurves_z$ (resp., $\fincurves^w)$ be the set
 of $\gamma \in \fincurves$ with $\gamma(0) = z$ (resp., $\gamma(t_\gamma) = w$). 
 We let $\fincurves_z^w = \fincurves_z \cap \fincurves^w$ and we 
define $\fincurves_z(D),\fincurves^w(D),\fincurves_z^w(D)$
similarly.

If $\gamma,\gamma^1 \in \fincurves$ with
$\gamma(t_\gamma) = \gamma^1(0)$, we define
the concatenation $\gamma \oplus \gamma^1$ by
$t_{\gamma \oplus \gamma^1} = t_\gamma +
t_{\gamma^1}$ and
\[     \gamma \oplus \gamma^1(t) =
      \left\{ \begin{array}{ll}
    \gamma(t), & 0 \leq t \leq t_\gamma \\
     \gamma^1(t - t_\gamma), & t_\gamma \leq
               t \leq t_\gamma + t_{\gamma^1}.
  \end{array}  \right.  \]
For every $w$, the map $(\gamma,\gamma^1) \mapsto
\gamma \oplus 
\gamma^1$ is continuous from $\fincurves^w \times
\fincurves_w$ to $\fincurves$.

  Suppose $f:D \rightarrow 
D'$ is a conformal transformation and $\gamma
\in \fincurves(D)$.  Let  
\[    s_t =s_{t,\gamma} =  
 \int_0^t  |f'(\gamma(s))|^2 \; ds . \]
If $s_{t} < \infty$ for all $t < t_\gamma$,
we define $f\circ \gamma$ by  
$f \circ \gamma(s_t) = f(\gamma(t))$.  If $s_{t_\gamma}
< \infty$ and $f$ extends continuously to the endpoints of $\gamma$,
 then $f \circ \gamma \in \fincurves(D')$
and $t_{f \circ \gamma} = s(t_\gamma)$.
If $\mu$ is a measure supported on the set of curves $\gamma$
 in $\fincurves(D)$ such that $f \circ \gamma$ is
well defined and in $\fincurves(D')$, 
 then $f \circ \mu$ will denote
the measure
\[     f \circ \mu(V) = \mu [\{\gamma:
    f \circ \gamma \in V \}] . \]

If $\gamma \in \fincurves$, we let $\gamma^R$ denote
the time reversal of $\gamma$, i.e., $t_{\gamma^R}
= t_\gamma$ and $\gamma^R(s) = \gamma(t_\gamma - s),
0 \leq s \leq t_\gamma$.  Similarly if $\mu$ is measure
on $\fincurves$, we define the measure $\mu^R$ in
the obvious way.

Suppose $\{\mu_D\}$ is a family of measures indexed
by a family of domains $D$
in $\C$.  We say that $\mu_D$ satisfies
the {\em restriction property} if
\begin{itemize}
\item $\mu_D$ is supported on $\fincurves(D)$;
\item if $D' \subset D$, then $\mu_{D'}$ is  
 $\mu_D$  restricted to the curves in $\fincurves(D')$. 
 \end{itemize} 
Note that if $\mu$ is any measure on $\fincurves$ and $\mu_D$ is
defined as $\mu$ restricted to $\fincurves(D)$, then 
the family $\{\mu_D\}$ satisfies the restriction property.  Conversely, suppose that 
\begin {itemize}
\item $\{\mu_D\}$ satisfies the restriction property
\item $D_n$ is an increasing sequence of domains whose union is $\C$
\item $\mu = \lim_{n \rightarrow \infty} \mu_{D_n}$.
\end {itemize}
Then, for each $D$, $\mu_D$ is $\mu$ restricted to $\fincurves(D)$.

If $A$ is any compact set, we define 
$\rad(A) = \sup\{|z| : z \in A\}$.  The
{\em half-plane capacity} of a subset $A$ of ${\overline \H}$
is defined by 
\begin{equation}  \label{jan23.5}
 \hcap(A) = \lim_{y \rightarrow \infty}
     y\,\E^{iy}[\Im(B_{\rho_A})] , 
\end{equation}
where $\rho_A = \inf\{t: B_t \in A \cup \R\}.$
It is not difficult to see that the limit exists, 
 satisfies the scaling rule $\hcap(rA)
= r^2\,   \hcap(A)$,  and  is monotone
in $A$.  If $A$ is such that $\Half \setminus
A$ is simply connected, then we use $\tilde g_A$
to denote the unique conformal transformation of
$\Half \setminus A$ onto $\Half$ such that
$\tilde g_A(z) - z = o(1)$ as $z \rightarrow \infty$.
Then $\tilde g_A$ has an expansion at infinity
\[  \tilde g_A(z) = z + \frac{\hcap(A)}{z} 
    + O(|z|^{-2}) . \]
 Since
$z \mapsto z + (1/z)$ maps $\{z \in \Half: |z| > 1\}$
conformally onto $\Half$, we can see that
$\hcap(\overline {\Disk_+}) = 1$.

If $\eta:[0,\infty) \rightarrow \C$ is a curve,
we will sometimes write $\tilde g_t$ for
$\tilde g_{\eta[0,t]}$ and define $g_t (z)
= \tilde g_t (z) - g_t (\eta_t)$ as in the introduction.

\section{Brownian bridges, Brownian bubbles}

We will start defining some bridge-type Brownian measures on curves that 
we will use. These are measures on Brownian paths with 
prescribed starting point and prescribed terminal
point. Since we are interested in conformally invariant 
properties, the standard bridges with prescribed time duration are not
well-suited.  

In our notation,
$\mu_D(z,w)$ will always be a measure on Brownian paths that remain in the 
domain $D$, that start at $z$ and end at $w$, but this notation
 will have different meanings depending
on whether $z,w$ are boundary or interior points of the domain $D$.  We
hope this will not cause confusion.
Since the content of this section is rather standard,
 we will just review these definitions.
The excursion measures have been defined in \cite {LW2,Vi,LSWrest}, 
the bubble measures in \cite {LSWrest}.

\subsection {First definitions}
\subsubsection {Interior to interior}
Let $\mu(z,\cdot;t)$ denote the law of 
a standard complex Brownian motion $(B_s, 0 \le s \le t)$,
with $B_0 = z$, viewed as an element of
$\fincurves$. We can write
\[
 \mu(z,\cdot;t) = \int_\C \mu(z,w;t) \; dA(w) , 
\]
where $A$ denotes area and $\mu(z,w;t)$ is a measure
supported on $\gamma \in \fincurves_z^w$ with
$t_\gamma = t$.  In other words, $\mu (z,w;t)$ is  $|\mu(z,w;t)|$ times the law 
$\mu^\normed (z,w;t)$ of the Brownian bridge from $z$ to $w$ in time $t$,
where $|\mu(z,w;t)| =
(2 \pi t)^{-1} \exp\{-|z-w|^2/(2t)\}$.

The measure $\mu(z,w)$ is defined by
\[  \mu(z,w) = \int_0^\infty \mu(z,w;t) \; dt . \]
This is a $\sigma$-finite measure (the integral explodes at infinity so that the 
total mass of large loops is infinite; when $z=w$, it also
diverges at $0$).

The measure $\mu (z,z)$ is an infinite measure on Brownian loops that 
start and end at $z$.  We can write
$$ 
\mu (z,z) = \int_0^\infty \frac {1}{2 \pi t} \, \mu^\normed (z,z;t)\; dt.
$$
where $\mu^\#(z,z;t)$ is the usual probability measure of a Brownian
 bridge from $z$ to $z$.

If $D$ is a domain and $z,w \in D$, we define
$\mu_D(z,w)$ to be $\mu(z,w)$ restricted to 
 $\fincurves(D)$.  For fixed $z, w$, the family 
$\{ \mu_D (z,w), D \supset \{ z, w \} \}$
clearly satisfies the
restriction property.

If $z\neq w$, and if the
domain $D$ is such that
a Brownian motion in $D$ eventually
exits $D$, then $|\mu_D(z,w)| < \infty$.
In fact, 
$$|\mu_D(z,w)| = \frac { G_D(z,w)}\pi ,$$
 where $G_D$ denotes the Green's function normalized so that
$G_\Disk(0,z) = -\log |z|$. 
 Note that $\mu_D(z,z)$ is well defined
and has infinite total mass.
The reversibility of the Brownian bridge immediately implies
 that $[\mu_D(z,w)]^R = \mu_D(w,z)$.  

\subsubsection {Interior to boundary}

Let $D$ be a connected domain in $\C$ whose
boundary is a finite union of  curves (we allow the curves
to be in the sphere and for infinity to
be a boundary point).  We will call
$\p D$ nice if it is piecewise analytic, i.e., if it is
a finite union of analytic curves.  A nice boundary point will be
any point at which the boundary is locally an 
analytic curve.

Let $B$ be a Brownian motion starting at $z \in D$
and stopped at its exit time of $D$,  i.e., at
 $$\tau_D = \inf\{t: B_t \not\in D\}.$$
Define $\mu_{D}(z,\p D)$ to be the law of
$(B_t, 0 \leq t \leq \tau_D)$.
If $D$ has a nice boundary we can write
\[  \mu_{D}(z,\p D) = \int_{\p D} \mu_D(z,w)
  \; |dw| , \]
where $\mu_D(z,w)$ for $z \in D$ and $w \in \p D$ denotes
a measure supported on $\fincurves_z^w(D)$ with
total mass $H_D(z,w)$, where $H_D(z,w)$ denotes
the usual Poisson kernel.  
The normalized probability measure $\mu_D^\normed(z,w)$ is the law of 
Brownian motion conditioned to exit $D$ ``at $w$''.

\subsection {First properties}

\subsubsection {Conformal invariance}

It is well known that planar Brownian motion is
conformally invariant.  In our interior to interior notation, this
can be phrased as follows. Suppose $f:D \rightarrow
D'$ is a conformal transformation and $z,w$ are two 
interior points in $D$.
Then,
\begin {equation}
f \circ \mu_D(z,w) = \mu_{f(D)}(f(z),f(w)) .
\end {equation}
If $z \neq w$, this is  a combination of the two
classical results:
$G_{f(D)}(f(z),f(w)) = G_D(z,w)$ and $[f \circ \mu_D]^\normed
(z,w) = \mu_{f(D)}^\normed(f(z),f(w)).$  For $z=w$ (in which
case the measures are infinite), one can prove this by taking a
 limit.

Similarly, in the interior to boundary case, if $z \in D$ is an interior 
point and $w$ a boundary point, and if both $w$ and $f(w)$ are nice, then 
\begin {equation}
f \circ \mu_D (z,w) = |f'(w)| \  \mu_{f(D)} ( f(z), f(w) ).
\end {equation}
This is a consequence of the two
relations: $H_D(z,w) = |f'(w)| \; H_{D'}(f(z),
f(w))$ and $f \circ \mu_D^\normed(z,w) =
\mu_{f(D)}^\normed(f(z),f(w))$.   
It implies  that one can define
the probability measure $\mu_D^\normed(z,w)$
for any simply connected $D$ and any boundary point (i.e. prime end) $w$
by conformal invariance. For instance, 
it suffices to put
 $\mu_D^\normed(z,w) = f \circ \mu_{\Disk}^\normed(0,1)$
where $f:\Disk \rightarrow D$ is the conformal
transformation with $f(0) = z$ and $f(1) = w$.  

\subsubsection {Regularity}

Note that the measures $\mu_D(z,w)$ are continuous functions of $z,w$
in the Prohorov metric. For instance, for any two interior points
$z_0 \not= w_0$ in the fixed domain
 $D$, the mapping $(z,w) \mapsto \mu_D (z,w)$ is continuous 
at $(z_0, w_0)$.
This can for instance be proved using a coupling argument.

Similarly, it is not difficult to show in the interior to boundary case
that for a fixed boundary point $w$, the mapping 
$z \mapsto \mu_D (z, w)$ is continuous.
When one wishes to let $w$ vary, one can for instance first note that 
$w \mapsto \mu_{\Disk} (0, w)$
is clearly continuous on the unit circle. Furthermore, for a conformal map
$f$ from $\Disk$ onto $D$, the derivative $f'$ is uniformly bounded when restricted to any $r \Disk$
for $r <1$, so that one can control the variation of the time-parametrization. We will discuss this in more 
detail later in the (slightly more complicated) case of the excursion measures.

\subsubsection {Relation between the two}

If $z,w$ are distinct points
in  $D$, then the normalized interior to interior
measure $\mu^\#_D(z,w)$ can be given as a limit of boundary measures.
Let $D_\epsilon = \{z' \in D: |z' - w| > \epsilon \}$,
and let $\nu_\epsilon$ denote $\mu(z,\p D_\epsilon)$ restricted
to curves whose terminal point is distance $\epsilon$ from
$w$.  As $\epsilon \rightarrow 0+$, $|\nu_\epsilon| \sim
G_D(z,w) \, [\log(1/\epsilon)]^{-1}  $ and
  $\nu_\epsilon^\# \rightarrow
\mu^\#_D(z,w)$.

The interior to boundary measure can also be viewed as the limit of an appropriately rescaled 
interior to interior measure:
If $w_n \in D$ and $w_n \rightarrow w$ where $w \in \p D$, 
then it is not hard to show
that the corresponding probability measures
converge $\mu^\normed_D(z,w_n) \rightarrow
\mu^\normed_D(z,w)$,  for instance using a coupling
argument.
Also, if $w$ is a nice boundary point,
and ${\bf n}_w$ denotes
the inward normal at $w$, then as $\epsilon
\rightarrow 0+$, 
$$
G_D(z,w - \epsilon {\bf n}_w) \sim
2 \pi  \epsilon 
H_D(z,w)
$$  
(the multiplicative constant can be worked out immediately using
 the case $D = \Disk, z=0,w=1$).
Hence,
\begin {equation}  \label{apr7.1}
  \lim_{\epsilon \rightarrow 0+}
        \frac{1}{2  \epsilon} \, \mu_D(z,
    w - \epsilon {\bf n}_\epsilon) =
       \mu_D(z,w),
\end {equation}
for any interior point $z$ and any nice boundary point of $w$.

\subsection {Excursion measures}

\subsubsection {Definition and conformal invariance}
Suppose   that $D$ is a nice domain, and
that $z$ and $w$ are different nice boundary points
of   $D$. We will define the Brownian measure on paths from
$z$ to $w$ in $D$. This Brownian excursion measure
$\mu_D (z,w)$  can be defined 
by various means (see e.g. \cite {LW2,Vi,LSWrest}).
It can be viewed as limits of the previous measures:
\begin {equation}
\label {c1}
  \mu_D(z,w) = \lim_{\epsilon \rightarrow 0+}
  \frac{1}{2  \epsilon^2} \mu_D(z + \epsilon {\bf n}_z,
   w + \epsilon{\bf n}_w) = \lim_{\epsilon \rightarrow 0+}
      \frac 1 { \epsilon} \mu_D(z+ \epsilon {\bf n}_z, w) . 
\end {equation}
Again we can write $\mu_D(z,w) = H_D(z,w) \, \mu^\normed_D(z,w)$ where 
\[
 H_D(z,w) = \lim_{\epsilon
\rightarrow 0+} \epsilon^{-1} H_D(z + \epsilon {\bf n}_z,w). 
\]
Under this normalization $H_\Half(0,x) = 1/ (\pi x^2)$.

The
probability measures $\mu^\normed_D(z,w)$ are conformally invariant, i.e., 
$$
f \circ \mu^\normed_D (z,w) = \mu^\normed_{f(D)} (f(z), f(w))
$$ 
for a conformal transformation such that the four boundary points $z$, $w$,
$f(z)$ and $f(w)$ are nice. This shows that one can define 
$\mu^\normed_D(z,w)$ by conformal invariance even
if $z,w$ are not nice boundary points.

It is sometimes
easier to consider $\mu_\Half^\normed(0,\infty)$ where $\Half$ denotes the upper half-plane.
 This is the distribution of $\Half$-excursions,  which are Brownian motions in 
 the first component and independent three-dimensional Bessel processes in 
 the second component,  see \cite {Vi, LSWrest}. 
  One could choose this as 
 the  definition of $\mu_\Half^\normed (0, \infty)$,  define the measures
 $\mu_D^\normed (z,w)$ by conformal invariance,  and define the 
measures $\mu_D (z,w)$ by multiplying by the total 
 mass, and finally verify (\ref{c1}).   (The measure $\mu_\Half^\normed(0,\infty)$ 
 is not supported on $\fincurves$ since curves under this measure have infinite time duration;
 however, this does not present a problem.  In particular, the   image of
$\mu_\Half^\normed(0,\infty)$ under a conformal transformation onto a bounded domain is
supported on paths of finite time duration.)

 If $f: D \rightarrow D'$ is a conformal transformation, and $z,w,f(z),f(w)$ are nice boundary points, then
\cite {Vi,LSWrest}
\begin {equation}
     f \circ \mu_D(z,w) = |f'(z)| \; |f'(w) | \;
    \mu_{f(D)}(f(z),f(w))
\end {equation}
The ``integrated measure''
\[    \mu_{\p D} := \int_{\p D} \int_{\p D}
     \mu_D(z,w) \; |dz| \; |dw| , \]
is therefore conformally invariant: 
\begin {equation}
f \circ \mu_{\p D}
= \mu_{\p f(D)}
\end {equation}
as was pointed out in \cite {LW2}.

\subsubsection {Regularity}

We now study the regularity of the excursion measures with respect to the domain $D$.
For this we will need some simple lemmas. 

\begin {lemma}
\label {sl1}
For any simply connected
domain $D$ and any two distinct points $w$ and $w'$
on the boundary of $D$, the expected time spent in an open subset $U$ of $D$
 by an excursion 
defined under the probability measure $\mu_D^\normed (w,w')$ 
is bounded from above by $2 \, {\rm area} (U) / \pi$.
\end {lemma}

\Pf
If $z \in D$, let $G^\normed_{D}(w,w';z)$
denote the Green's function for $\mu^\normed_D(w,w')$.
This can be obtained as the limit of
$G^\normed_{D}(w_n,w';z)$ where $w_n$ is a
sequence of points in $D$ converging to $w$.
If $f:D \rightarrow
D'$ is a conformal transformation,  then
\[    G^\normed_{D}(w,w';z) =
   G^\normed_{f(D)}(f(w),f(w');f(z)) . \]
Also 
\begin{equation}  \label{jan19.1}
 G^\normed_\Half(0,\infty;z) = 
      \lim_{\epsilon \rightarrow 0+} \frac{\Im(z)}
              {2 \, \epsilon}
       \log \frac{\Re(z)^2 + (\epsilon + \Im(z))^2}
     {\Re(z)^2 + (\epsilon - \Im(z))^2 } =
      2 \; \frac{\Im(z)^2}{|z|^2} \leq   2   . 
\end{equation}
By conformal invariance, we get
$G^\normed_D(w,w';z) \leq 2 $ for all simply
connected $D$ and all $w,w',z$. This readily implies the lemma.
\endpf

\begin{lemma} 
There is a constant $c < \infty$ such that the following
holds. Suppose $D,D'$ are simply connected
domains and $f:D \rightarrow D'$ is a conformal
transformation with $|f(z_1)-z_1| \leq \delta  \leq 1$
 for all $z_1 \in D$.  Then for any path $\gamma$ in $D$, 
\begin {eqnarray*}
 \finmetric(\gamma,f\circ \gamma)  
 &\leq&
 c \; [ \; \delta + \delta^{1/2} t_\gamma +
  \int_0^{t_{\gamma}} 1\{\dist(\gamma(s),
\p D) \leq \delta^{1/2} \}  \; ds   \\
&& + 
   \int_0^{t_{f \circ \gamma}}
      1\{\dist(f \circ \gamma(u),
   \p D') \leq  c \delta^{1/2} \} \; du\; ] . 
   \end {eqnarray*}
\end{lemma}

\Pf For any $\gamma$ let
\[  \theta_\gamma(s) = \int_0^{s}
    |f'(\gamma(r))|^2 \; dr . \]
Then,
\[ \finmetric(\gamma,f\circ \gamma) 
   \leq \sup_{0 \leq s \leq t_\gamma}
          |s- \theta_\gamma(s)|
    + \sup_{0 \leq s \leq t_\gamma}
             |\gamma(s) - f(\gamma(s))| . \]
The second term on the right is bounded above by
$\delta$ and the first term is bounded above
by
$  \int_0^{t_\gamma}  Y_s \; ds , $
where $Y_s = | \, |f'(\gamma(s))|^2 - 1 \, |$.
Let 
$$\tilde D = \tilde D_\delta = \{
z \in D: \dist(z,\p D) \leq \delta^{1/2} \}.$$
For $z \in D \setminus \tilde D$, a standard estimate
gives $|f'(z) - 1| \leq c\, \delta/\dist(z,\p D)
\leq c \, \delta^{1/2}$.  Hence,
\[    \int_0^{t_\gamma} Y_s \, 1_{\gamma(s)
     \not\in \tilde D} \, ds \  \leq c \,
    \delta^{1/2} \,  t_\gamma 
    . \]
For the other part, write
\[  \int_0^{t_\gamma} Y_s \, 1_{\gamma(s)
    \in \tilde D} \, ds
  \leq \int_0^{t_\gamma} 1_{\gamma(s)
    \in \tilde D} \, ds  +
    \int_0^{t_\gamma} |f'(\gamma(s))|^2 \;
    1_{\gamma(s) \in \tilde D} \, ds . \]
The first term on the right hand side is the
amount of time that $\gamma$ spends within
distance $\delta^{1/2}$ of the boundary. 
Since $|f(z) - z| \leq \delta$, the second term
on the right is less than the amount of time
that $f \circ \gamma$ spends within distance
$2\delta^{1/2}$ of $\p D'$.  
Combining these estimates gives the lemma.
\endpf. 
 
\begin {lemma}
\label {sl3}
Suppose that $D \subset \Disk_+$ is a simply connected domain with
$\Disk_+ \setminus D \subset \delta \Disk_+$ for some $\delta > 0$.
Let  $z,z',w $  on $\partial D$ 
with $|z|= |z'|= 1$, $|w| \le \delta$ and $|z-z'| \le \delta$.
Then, the distance between 
$\mu_D^\normed (z,w)$ and $\mu_{\Disk_+}^\normed (z', 0)$ goes to zero with $\delta$,
uniformly with respect to the choice of $w,z,z'$ and $D$.
\end {lemma}

\Pf
Let $f$ denote the conformal mapping from $D$ onto $\Disk_+$ such that $f(w)=0$, $f(z)=z'$ 
and $|f'(i)| = 1$. It is standard that for some constant $c$,
$| f(x) - x|\le c \delta$ for all $x \in D$.
The total area in $D$ or $D'$ of the set of points that is at distance less than 
$\delta^{1/2}$ from the boundary is no larger than $c' \delta^{1/2}$. Hence, a combination 
of the two previous lemmas shows that
$$
 \E[ \finmetric(\gamma,f\circ \gamma)] \leq
          c \, \delta^{1/2} , \]
where the expectation is with respect
to $\mu_D^\#(z,w)$.  Since the law of $f \circ \gamma$ is 
$\mu_{\Disk_+}^\normed (z',0)$, the lemma follows.
\endpf

\subsection{Brownian Bubbles}
 
\subsubsection {Definition}

The Brownian bubble measure in $\H$ at the origin is the $\sigma$-finite measure
\begin{equation}  \label{jan22.1}
    \bubble_\Half(0) = \lim_{z \rightarrow
    0} \frac{\pi}{\Im(z)}\, \mu_\Half(z,0)
   \;\;\;\; (z \in \Half)
\end{equation}
or equivalently (see (\ref{apr7.1})),
\begin{equation} \label{jan23.4}
 \bubble_\Half(0) = \lim_{z,w \rightarrow
0} \frac{\pi}{2\, \Im(z)\, \Im(w)} \, \mu_\Half(z,w)
\;\;\;\;  (z,w \in \Half). 
\end{equation}
When we speak of the limit, we mean that for every
$r > 0$, if we restrict the measures 
on the right to loops
that intersect the circle of radius $r$ (so that this is a finite 
measure), then
the limit exists and equals $\bubble_\Half(0;r)$
which is $\bubble_\Half(0)$ restricted to loops
that intersect $\{|z| = r\}$. 
It is not hard to show the limit exists
and the normalization is chosen so that
$| \bubble_\Half(0;r)| = 1/r^2$.
If $r > 0$ and $f_r(z)
= rz$, then $\bubble_\Half(0)$ satisfies
the scaling rule
\[     f_r \circ \bubble_\Half(0) =
             r^2 \; \bubble_\Half(0) . \]
We can also define $\bubble_D(z)$ for other domains,
at least if $\p D$ is smooth near $z$, 
using conformal covariance,
\[    f \circ \bubble_D(z) = |f'(z)|^2 \;
       \bubble_{f(D)}(f(z)) . \]
These measures satisfy the restriction property:
if $D' \subset D$, then $\bubble_{D'}(0)$ is
$\bubble_D(0)$ restricted to loops that are in
$D'$.

Suppose $D \subset \Half$ is a simply connected
domain containing $r \Disk_+$ for some
$r > 0$, and let $A$ be the image of $\Half \setminus
D$ under the map $z \mapsto - 1/z$.  
Then \cite[(7.2)] {LSWrest} tells us that
the $\bubble_\Half(0)$ measure of the set of loops
that do not stay in $D$ is $\hcap(A)$.
 The reader can check
 that both the definition  in the present paper and
the definition in \cite{LSWrest} give measure $1$ to
the set of loops that intersect the unit circle, and hence
 the two 
definitions use the same normalization. 
In particular, this shows immediately that 
\begin {equation}
\label {schw}
\bubble_\Half (0) [\{ \gamma \ : \gamma (0, t_\gamma) \not\subset D \} ] 
= \frac {- S_\Phi (0)}6
,\end {equation}
where $\Phi$ is a conformal map from $D$ onto $\H$ that 
keeps the origin fixed, say, and $S_\Phi$ denotes the Schwarzian derivative
\[   S_\Phi(z) = \frac{\Phi'''(z)}{\Phi'(z)} - \frac{3 \Phi''(z)^2} {2 \Phi'(z)^2} . \]

\subsubsection {Path decomposition}

The next proposition relates $\bubble_\Half(0)$
to excursion measures. This expression for
$\bubble_\Half(0)$ splits the bubble 
at  the point
$s e^{i \theta}$ at which its distance to the origin is maximal.

\begin{proposition} 
One has
\begin{eqnarray}
  \label{jan23.1}
  \bubble_\Half(0) &=& \pi \int_0^\infty \int_0^\pi
    [\mu_{r \Disk_+}(0, r e^{i \theta})
    \oplus \mu_{r \Disk_+}(re^{i \theta},0)]\; r \; d\theta \; dr\\
   \label {bis}
    &=&
     \int_0^\infty 
  \frac {4}{\pi r^3} \int_0^\pi
    [\mu_{r \Disk_+}^\#(0, r e^{i \theta})
    \oplus \mu_{r \Disk_+}^\#(re^{i \theta},0)]\; 
 \sin^2 \theta  \; d\theta \;   dr. 
\end{eqnarray}
\end{proposition}

\Pf Let $r > 0, \delta > 0$.
By the strong Markov property,
\begin{eqnarray*}
 \lefteqn {\bubble_\Half(0;r) - \bubble_\Half(0;r+ \delta)} \\  
 & = & \lim_{\epsilon \rightarrow 0+}
     \frac \pi \epsilon
      \int_0^\pi
              [\mu_{r \Disk_+}(\epsilon i, r e^{i \theta})
    \oplus \mu_{(r+\delta) \Disk_+}(re^{i \theta},0)]
 \;  r \; d\theta .  \\
   & = & \pi   \int_0^\pi
[\mu_{r \Disk_+}(0, r e^{i \theta})
    \oplus \mu_{(r+\delta) \Disk_+}(re^{i \theta},0)]
  \; r \;  d\theta .
\end{eqnarray*}
But as 
$$\lim_{\delta \to 0+}
\delta^{-1}  \mu_{(r+\delta) \Disk_+}(re^{i \theta},0) = \mu_{r \Disk_+}(re^{i \theta},0),
$$
we get that 
\[   \frac{d}{dr}  \bubble_\Half(0;r) 
       = - \pi \int_0^\pi [\mu_{r \Disk_+}(0, r e^{i \theta})
    \oplus \mu_{r \Disk_+}(re^{i \theta},0)]\; r \; d\theta,\] which gives (\ref{jan23.1}). 
Identity  (\ref {bis}) follows from the fact (see Lemma \ref{apr3.lemma1}) that 
$$|\mu_{r \Disk_+}(0, r e^{i \theta})| = r^{-2} \;
  |\mu_{\Disk_+}(0,e^{i \theta})| = 
  \frac {2 \sin \theta}{\pi r^2}. 
$$  
\endpf

Note that 
\[  \frac d {dr} | \bubble_\Half(0;r)|
      = - 4/(\pi r^3)  \int_0^\pi \sin^2 \theta \; d\theta = -2 / r^3 , \]
which is consistent with $| \bubble_\Half(0;r)|
= r^{-2}$.

Similarly, one can decompose the Brownian bubble measure 
in $\H$ at the point where the imaginary part is maximal.
This gives a joint description of the real and imaginary parts using 
one-dimensional excursions and Brownian bridges (as briefly 
mentioned in \cite {LSWrest}).

\section{(Unrooted) loop measure}

\subsection {Definition, restriction and conformal invariance}
We will now define the most important object for
this paper, the Brownian loop measure $\loopmeasure$.
Let $\percurves$ be the set of loops, i.e., the
set of $\gamma \in \fincurves$ with $\gamma(0) = \gamma(t_\gamma)$. 
Such a $\gamma$ can also be considered as a  function with domain
 $(-\infty,\infty)$ satisfying $\gamma(s) = \gamma(s + t_\gamma)$.

Define $\theta_r: \percurves
\rightarrow \percurves$ by $t_{\theta_r\gamma}
= t_\gamma$ and $\theta_r\gamma(s) = \gamma(s+ r)$.
We say that two loops $\gamma$ and $\gamma'$ are equivalent if  
 for some $r$, $\gamma' = \theta_r \gamma$.
We write $[\gamma]$ for the equivalence class of $\gamma$. 
Let $\unrootedcurves$ be the set of {\em unrooted
loops}, i.e., the equivalence classes in $\percurves$.
Note that $\unrootedcurves$ is a metric
space under the metric
\[ \unrootedmetric(\gamma,\gamma') = \inf_{r \in [0, t_\gamma]}
      \finmetric(\theta_r\gamma, \gamma'). \]

Any measure supported on $\percurves$ gives a measure
on $\unrootedcurves$ by ``forgetting the root'', i.e.,
by considering the map $\gamma \mapsto [\gamma]$.
If $D$ is a domain, we define $\percurves(D),\unrootedcurves(D)$ to be the set of loops that lie 
entirely in $D$, i.e., $\gamma[0,t_\gamma] \subset D$.

We define the{ \em Brownian loop measure} $\loopmeasure$ on $\unrootedcurves$ by
\begin {equation}
\loopmeasure 
= \int_\C \frac{1}{t_\gamma} \, \mu (z,z)\; {dA(z)}=
\int_\C \int_0^\infty   \frac {1}{2 \pi t^2 }\mu^\normed (z,z; t)\; dt \; dA (z),
\end {equation}
where $dA$ denote the Lebesgue measure on $\C$.
We insist on the fact that the measure $\loopmeasure$
 is a measure on {\em unrooted} loops.  
 
We will call a Borel measurable
function $T: \percurves \rightarrow
[0,\infty)$ a {\em unit weight} if for
every $\gamma \in \percurves$,
\[   \int_0^{t_\gamma} T(\theta_r\gamma) \; dr
    = 1 . \]
One example of a unit weight is $T(\gamma) =
1/t_\gamma$.  Note that   $\loopmeasure$ satisfies
\begin{equation}  \label{jan14.1}
  \loopmeasure = \int_\C T\; \mu(z,z) \;
    dA(z)  
\end{equation}
(considered as a measure on $\unrootedcurves$) for any  unit
weight $T$. 

If $D$ is a domain,  we define $\loopmeasure_D$
to be  $\loopmeasure$ restricted to the curves 
in $\unrootedcurves(D)$; this is the same as the
right-hand side of (\ref{jan14.1}) with
$D$ replacing $\C$ and 
$\mu_D(z,z)$ replacing $\mu(z,z)$.
  By construction, the family
$\{\loopmeasure_D\}$ satisfies the restriction
property.
Not as obviously, these measures are also
conformally invariant:

\begin {proposition}
If $f: D
\rightarrow D'$ is a conformal transformation,
then $f \circ \loopmeasure_D = \loopmeasure_{f(D)}$.
\end {proposition}

\Pf
Showing this requires two observations.  One,
which we have already noted, is the conformal
invariance of interior to interior measures,
$f \circ \mu_D(z,z)
= \mu_{f(D)}(f(z),f(z))$.  The other is the fact that
 we can define  a 
unit weight $T_f$ by $T_f(\gamma) = 1/t_\gamma$
if $\gamma \not\in \percurves(D)$, and if $T
\in \percurves(D)$, $T_f(\gamma) = |f'(\gamma(0))|^2/
t_{f \circ \gamma}$. To check that this is a
unit weight, note that
\[  \int_0^{t_\gamma} T_f(\theta_r\gamma) \; dr
  = (1/t_{f \circ \gamma}) \int_0^{t_\gamma}
          |f'(\gamma(r))|^2 \; dr = 1 . \]
Therefore,
\begin{eqnarray*}
f \circ \loopmeasure_D & = &
   f \circ \int_D T_f \; \mu_D(z,z) \, dA(z)\\
  & = &  \int_D
(1/t_{f \circ \gamma}) \; 
     |f'(z)|^2\; f \circ \mu_D(z,z) \; dA(z) \\
  & = &   \int_D 
   (1/t_{f \circ \gamma}) \;
      \mu_{D'}(f(z),f(z)) \; [|f'(z)|^2\;dA(z)] \\
   & = & \int_{D'} T \, \mu_{D'}(w,w)\; dA(w)
  = \loopmeasure_{D'}.
\end{eqnarray*}
Here $T$ denotes the simple unit weight
$T(\gamma) = 1/t_\gamma$.
\endpf

Note that the same argument shows that $\loopmeasure$ is invariant under the inversions 
$z \mapsto 1 / (z-z_0)$ for all fixed $z_0$.

\subsection {Decompositions}

The definition of $\loopmeasure$ makes it conformally
invariant and hence independent of the choice of
coordinate axes.  It will be however convenient to have expressions for
$\loopmeasure$ that do depend on the axes. 
We will  write the measure on unrooted loops $[\gamma]$ as a measure on rooted
 loops by choosing the representative $\gamma$ whose initial point is the (unique) point on the
loop of minimal imaginary part (the same works of course also for the maximal imaginary part).
Note that this choice of ``root'' of the loop is not conformally invariant. 

\begin{proposition}  \label{prop.jan23.1}
$$
  \loopmeasure  =    \frac 1 {2\pi} \; \int_{-\infty}^\infty
    \int_{-\infty}^\infty \bubble_{\Half + iy}
    (x+iy) \; dx \; dy 
$$
\end{proposition}

\Pf There are various simple ways to prove this. The 
main point is to get the multiplicative constants right. We therefore opt for 
a self-contained elementary proof that does not rely on other multiplicative 
conventions (i.e., excursions). 
 We start by recalling some facts about one dimensional Brownian motion.  
 Suppose $Y_t$ is a one-dimensional Brownian motion started at the origin. 
 Let $t^*$ be the time in $[0,1]$ at which $Y_t$ is minimal, let $M = Y_{t^*}$,
and let $\Psi = (Y_0 - M) \, (Y_1 - M)$. It is easy to see that the law of 
$t^*$ is the arcsine law with density $1/ ( \pi \sqrt { t (1-t)})$ on $[0,1]$.
Given $t^*$, $Y_0 - M$ and $Y_1-M$ are independent random variables with
the distribution of Brownian motion ``conditioned to stay positive''. 
It is not difficult to show that $\E[Y_0 - M \mid t^* = t] =
\sqrt{\pi t/2}$ and hence 
 that $\E [ \Psi] = 1/2$.

We now define a unit weight $T_\epsilon$ on
$\gamma \in \percurves$ that will approximate the Dirac
 mass at the time of the 
minimal imaginary part of $\gamma$.
  If $t_\gamma <
\epsilon$, then $T_\epsilon(\gamma)
= 1/t_\gamma$.   Suppose  $t_\gamma \geq
\epsilon$ and there is
a unique $r_0 \in [0,t_\gamma)$ such that
$\Im[\gamma(r)] < \Im[\gamma(t)]$ for
$t \in [0,t_\gamma) \setminus \{r_0\}$.  
Then $T_\epsilon(\theta_r \gamma) =
1/\epsilon$ for $r_0 - \epsilon \leq
r \leq r_0$ and $T_\epsilon(\theta_r \gamma)
=0$ for other $r_0 < t < r_0 + \epsilon$
(here $\gamma$ is considered as a periodic
function of period $t_\gamma$).  If no
such unique $r_0$ exists, set
$T_\epsilon (\gamma) = 1/t_\gamma$ (the
choice here is irrelevant since this is a
set of loops of measure zero).  Note that
the measures $\mu(z,z)$ are supported on 
loops for which a unique $r_0$ exists. 
It is easy to see that $T_\epsilon$ is a unit
weight, and hence for every $\epsilon$,
\[    \loopmeasure = \int_\C T_\epsilon
  \; \mu(z,z) \; dA(z) =  \lim_{\epsilon
\rightarrow 0+} 
   \int_\C \epsilon^{-1}
  \; \mu(z,z; \ge \epsilon) \; dA(z) , \]
where $\mu (z,z; \ge \epsilon )$ denotes $\mu(z,z)$
restricted to curves $\gamma$  with 
$t_\gamma  \geq  \epsilon$  and 
\[     \inf\{\Im(\gamma(t)); 0 \leq t
   \leq \epsilon\} = \inf\{\Im(\gamma(t)):
    0 \leq t < t_\gamma \} . \]
For fixed $\epsilon$, $\int_\C \epsilon^{-1}
  \; \mu(z,z;\geq \epsilon) \; dA(z)$ is the
same as $\loopmeasure$ restricted to curves
with $t_\gamma \geq \epsilon$.
Let us consider the measure $\epsilon^{-1} \; \mu(z,z;\geq \epsilon)$. For ease let $z=0$. 
Start a Brownian motion $B_t$ at $0$ and let it
run until time $\epsilon$; let us write $B_\epsilon
= \sqrt{\epsilon}w$.  We let $ - b
\sqrt{\epsilon} = \min\{\Im(B_t):
0 \leq t \leq \epsilon\}. $   Then given $B_t,
0 \leq t \leq \epsilon$, the remainder of the curve
is obtained from the measure $\epsilon^{-1}
\; \mu_{\Half - i b \sqrt {\epsilon}}
(0,w \sqrt{\epsilon})$.  As $\epsilon \rightarrow
0+$, this looks like $\epsilon^{-1} \, \mu_\Half
(i b \sqrt \epsilon, i(b+w) \sqrt{\epsilon})$,
which in turn has the same limit as
$\epsilon^{-1} \, b [b + \Im(w)] \, \mu_\Half
(i \sqrt \epsilon, i \sqrt \epsilon) .$
Hence (see (\ref{jan23.4})),
\begin{eqnarray*}
  \lim_{\epsilon \rightarrow 0+} 
     \epsilon^{-1} \;
\mu(0,0;\epsilon) & =  & \lim_{\epsilon
  \rightarrow 0+}  \E[b(b+\Im(w))] \;
  \epsilon^{-1} \mu_\Half
(i \sqrt \epsilon, i \sqrt \epsilon) \\
  & = & \frac {1}{2 \pi}
  \bubble_\Half(0) . 
\end{eqnarray*}
\endpf

The next proposition is similar. It gives an
expression for $\loopmeasure_\Half$ by associating
to an unrooted loop the rooted loop whose root
has maximal
absolute value.

\begin{proposition}  \label{jan23.prop2}
\[  \loopmeasure_\Half =
  \frac 1 {2 \pi}  \int_0^\infty \int_0^\pi  \bubble_{r \Disk_+}
    (r e^{i \theta}) \; d \theta \; r \, dr . \]
\end{proposition}

\Pf Let 
\begin {eqnarray*}
 \rect & =& \{x+iy: -\infty < x < \infty , 0 < y < \pi \} \\
 \rect_b & = & \{z \in \rect: \Re(z) < b\} 
 \end {eqnarray*}
 and let $\phi (z) = e^z$.  Conformal invariance tells us that
$\phi \circ \loopmeasure_\rect = \loopmeasure_\Half$.
But Proposition \ref{prop.jan23.1} (rotated ninety degrees)
and restriction tell us that
\[  \loopmeasure_\rect = \frac 1 {2 \pi}  \,  \int_{-\infty}^\infty
    \int_0^\pi \bubble_{\rect_x}(x+iy) \; dy \; dx. \]
The scaling rule for $\bubble$ gives
$\phi \circ  \bubble_{\rect_x}(x+iy) = e^{2x}
 \; \bubble_{e^x \Disk_+}(e^{x+iy})$.   Therefore
\begin{eqnarray*}
 \loopmeasure_\Half = \phi \circ \loopmeasure_\rect 
    & = & \frac 1 {2\pi}  \,
 \int_{-\infty}^\infty  \int_0^\pi 
            e^{2x}
 \; \bubble_{e^x \Disk_+}(e^{x+iy}) \; dy \; dx \\
  & = & \frac 1 {2\pi} \,
\int_0^\infty \int_0^\pi 
    \bubble_{r \Disk_+}(re^{iy} ) \; dy \; r \, dr . 
\end{eqnarray*}
\endpf

\noindent
{\bf Remark.}
If we combine this description with the invariance of the unrooted loop measure
under the inversion $z \mapsto -1/z$, we get that, when $r \to 0$, 
the measure $\loopmeasure_\Half$ restricted to those loops 
that intersect $r\Disk_+$ is close to a multiple of the Brownian bubble measure in $\H$ 
at the origin.

 \section {Bubbles  and loops}

The goal of this section is to derive the relation between the Poissonian cloud of 
loops that intersect a given curve and the Poisson point process of bubbles that we briefly 
described in the introduction.
In order to prove this, we need a clean generalization of the 
previous remark  to shapes other 
 than disks, and to show that the convergence holds uniformly over all shapes.

\subsection {Some estimates}

We will need some standard estimates about the Poisson kernel
on rectangles and half-infinite rectangles, or, more precisely,
on the images of these domains under the exponential map.

\begin{lemma} \label{apr3.lemma1}
There exist a constant $c$ such that if
$r \in (0,1/2) $ and $\theta,\varphi \in (0,\pi)$,
\begin{equation}  \label{C}
   |H_{\Disk_+}(r e^{i\theta},e^{i \varphi})  
                     - \frac 2 \pi \, r\, \sin \theta \, \sin \varphi| \leq
             c \, r^2 \, \sin \theta \, \sin \varphi , 
\end{equation}
\begin{equation}  \label{A}
 |H_{\Half \setminus \overline{\Disk_+}}(r^{-1} e^{i \theta},
   e^{i\varphi}) - \frac 2 \pi \, r \, \sin \theta \, \sin \varphi | \leq
    c \, r^2 \sin \theta \, \sin \varphi.
\end{equation}
\end{lemma}

\Pf The map $f(z) = -z -(1/z)$ maps $\Disk_+$ onto $\Half$.  Hence
\begin{eqnarray*}    H_{\Disk_+}(r e^{i\theta},e^{i\varphi}) & =  & 
      |f'(e^{i\varphi})| \, H_\Half(f(re^{i\theta}),f(e^{i\varphi})  )\\
    & = & 2 \, \sin \varphi \, H_\Half(f(re^{i\theta}),f(e^{i\varphi})).
\end{eqnarray*}
But if $|z| \geq 5/2, |x| \leq 2$,
\[   H_\Disk(z,x') =   \frac{\Im(z)}{\pi \; [(\Re(z) - x')^2 + \Im(z)^2]}
                          = \frac{\Im(z)}{\pi \, |z|^2} \, [1 + O(\frac{1}{|z|})]  , \]
and 
\[   f(r e^{i\theta}) = \frac{1}{r} \, e^{i(\pi - \theta)} \, +O(r), \]
\[  \Im[f(r e^{i \theta})] =    \frac{\sin \theta}{r}  + \sin \theta \, O(r) . \]
This gives the first expression, and the second is obtained from the
first using the map $z \mapsto -1/z$.  \endpf
 
\begin{lemma}  There exists a constant $c$ such that if $e^{-s} \in (3/4,1)$,
$r \in (0,1/2)$, 
and $\theta,\varphi \in (0,\pi)$, then
\begin{equation}  \label{B}
|H_{ \Disk_{+,r}}(e^{-s + i \theta},re^{i\varphi}) -
    \frac{4}{\pi} \, \sinh s \, \sin \theta \, \sin \varphi| \leq
       c \, r \, s \, \sin \theta \, \sin \varphi , 
\end{equation}
where $\Disk_{+,r} = \{z \in \Disk_+: |z| > r \}.$
\end{lemma}

\Pf  Separation of variables gives an exact form for the Poisson
kernel on a rectangle, and the logarithm maps $\Disk_{+,r}$
onto a rectangle.  Doing this we see, in fact,
that    
\[   H_{ \Disk_{+,r}}(e^{-s + i \theta},re^{i\varphi})   
                          = \frac{4}{\pi \, r}  \sum_{n=1}^\infty   \sin(n\theta) \, \sin(n \varphi) \,
                                \sinh(ns)  \;  \frac{r^n}{1 + r^{2n}} , \]
from which the estimate comes easily.  \endpf

\subsection{Bubble measure and loop measure}  \label{estimatesec}

Suppose $V_n$ is a sequence of  sets in $\Half$
with $u_n = \rad(V_n) \rightarrow 0$ and such
that $\Half \setminus V_n$ is simply connected.
Let $m_n$ be  $\loopmeasure_\Half$ restricted to loops that
intersect both $V_n$ and the unit circle.
We set $ h_n = \hcap (V_n)$.

\begin {proposition}
\label {p.tech}
When $n \to \infty$, 
$$
m_n  =  \frac {h_n}2 \, \bubble_\Half (0,1) ( 1 + o(1)),
$$
where $o(1)$ is uniformly bounded by a function of $u_n$ that goes to zero with $u_n$.
\end {proposition}

Note that scaling implies the corresponding results for the measures
restricted to paths that intersect any given circle $r \partial \Disk$.
\medbreak
\noindent
\Pf
 We will use the decomposition of the 
measures according to the point at which the loop (or the bubble) has maximal
absolute value. 
Recall that
\[
m_n = \frac 1 {2 \pi}  \int_1^\infty \int_0^\pi  \bubble_{r \Disk_+}
    (r e^{i \theta}| V_n) \; d \theta \; r \, dr, 
\] where $\bubble_{r \Disk_+}
    (r e^{i \theta}| V_n)$ denotes $\bubble_{r \Disk_+}
    (r e^{i \theta})$ restricted to loops that intersect $V_n$.
    Recall also that 
\begin {equation}
\label {recall}
 \bubble_\Half(0,1) =
     \int_1^\infty 
 \frac {4}{\pi r^3}
  \int_0^\pi
    [\mu_{r \Disk_+}^\#(0, r e^{i \theta})
    \oplus \mu_{r \Disk_+}^\#(re^{i \theta},0)]\; 
 \sin^2 \theta  \; d\theta \;   dr. 
 \end {equation}
It is not difficult to show that
\[     
[\bubble_{r \Disk_+} (r e^{i \theta}|V_n)]^\# 
\rightarrow
      \mu_{r\Disk_+}^\#(r e^{i \theta },0) \oplus
      \mu_{r\Disk_+}^\#(0,r e^{i \theta }), \]
uniformly on $\{1 \leq r \leq R\}$ and $\theta \in (0, \pi)$.
(Note that there is a conformal transformation 
$g: \Half \setminus  V_n
\rightarrow \Half$ with $\max|g(z) -z| \leq
c \,u_n $).  

We now focus on the total masses.
We claim that
\begin{equation}  \label{jan23.2}
 |\bubble_{ \Disk_+} 
    ( e^{i \theta}| V_n) |= 4\,
     h_n \sin^2 \theta \; [1+ O(u_n)]  . 
\end{equation}
By the scaling rules for $\bubble$ and $\hcap$ this
implies that for all $r \geq 1$,
\[ |\bubble_{ r\Disk_+} 
    ( re^{i \theta}| V_n) |= 4\;r^{-4} \;
     h_n \sin^2 \theta \; [1+ O(u_n)],\]
and the proposition follows, using (\ref {recall}).

To prove (\ref{jan23.2}), we first note that
$$
|\bubble_{ \Disk_+} 
( e^{i \theta}| V_n) |
 = \lim_{ \epsilon \to 0+} \frac {\pi}{\eps}
 \mu_{\Disk_+} ( \exp (- \eps + i \theta), \exp ( i \theta)) [ B \hbox { hits } V_n ].
$$
The estimates (\ref{A}) and (\ref{B})
 show that the following two measures are very close (for all large $R$, 
and $r$ is small):      
\begin {itemize}
\item
The measure of $1_{\sigma_u < T} \arg (B_{\sigma_u})$ when $B$ is defined under the 
measure $R \P^{iR}$.
\item
The measure of $1_{\sigma_u < T, \sigma_u < \sigma_1} \arg (B_{\sigma_u})$
 when $B$ is defined under the measure 
$$(2 \,
\sin \theta \sinh(\eps))^{-1} \P^{(1-\eps) \exp (i \theta)}.$$
\end {itemize}
After these hitting times, it is possible to ``couple'' the two paths up to their first 
hitting of $V_n$. After the hitting of $V_n$, we want to estimate the probability that 
the path go back to the unit circle without hitting $\R$ and that they hit it in the 
neighborhood of $\exp (i \theta)$. By (\ref{C}), this
will occur with a probability 
$$\frac {2}{\pi} \Im (B_{\rho_{V_n}}) \sin \theta d\theta .$$
Hence, we get finally that (recall that the estimates are uniform in $\theta$, $\epsilon$ and $R$)
\begin {eqnarray*}
|\bubble_{ \Disk_+}  ( e^{i \theta} \mid V_n ) |
& \sim & \lim_{\eps \to 0, R \to \infty} 
\frac {2 \, \pi}{\eps}  {R \sin \theta}{\sinh \eps}
\E^{iR} [ \frac {2}{\pi} \Im ( B_{\rho_{V_n}} ) \sin \theta ] \\
&\sim &
4 \sin^2 \theta \lim_{R \to \infty} R \E^{iR} [ \Im ( B_{\rho_{V_n}} )] 
\\
& \sim& 
4 h_n  \sin^2 \theta
\end {eqnarray*}
when $u_n \to 0$.
\endpf

\subsection{Bubble soup and loop soup}  \label{soupsec}

We define a {\em bubble soup} with
intensity $\lambda \geq 0$ to be a Poisson point process with intensity 
$\lambda \bubble_\Half$. One can also view it as a Poissonian 
sample from the measure 
$\lambda \bubble_\Half(0) \times ({\rm length})$ on $\fincurves_0^0(\Half)
\times [0,\infty)$. We can write a realization
of the bubble soup as a countable collection
$ {\cal U} = \{(\gamma_j,s_j)\}$.
Recall that the law of ${\cal U}$ is characterized by the fact that:
\begin {itemize}
\item
For any two disjoint measurable subsets $U_1$ and $U_2$ of
$\fincurves_0^0(\Half) \times [0,\infty)$, ${\cal U} \cap U_1$ 
and ${\cal U} \cap U_2$ are independent.
\item
The law of the number of elements in ${\cal U} \cap U$ is the Poisson law with 
mean $\lambda \bubble_\Half(0) \times ({\rm length}) [ U] $
(when this quantity is finite).
\end {itemize}
We will think of the bubble $\gamma_j$ as
being created at time $s_j$.
Clearly, with probability one $s_j \neq s_k$ for $j\neq k$.

A {\em Brownian loop soup} with intensity $\lambda$
is a Poissonian sample from the measure $\lambda \loopmeasure$.
We will use ${\cal L}_\C$ to denote a realization of the loop soup. 
A sample of the Brownian loop soup is a countable collection 
of (unrooted) Brownian loops in the plane. We will use ${\cal L}$ to denote 
the family of loops in ${\cal L}_\C$ that are in $\H$. This is the Brownian loop soup in the half-plane.
  
If $D \subset \Half$ is a domain, then we write 
\begin {itemize}
\item ${\cal L}(D)$ for the family of loops in ${\cal L}$ that are in $D$
\item
 ${\cal L}^\perp(D)$ for ${\cal L} \setminus {\cal L}(D)$, i.e.,
 the family of loops that intersect $\Half \setminus D$.
\end {itemize}
By definition, for any fixed $D$, the two random families 
${\cal L}(D)$ and ${\cal L}^\perp(D)$ are independent.

Note that the (law of the) families ${\cal L} (D)$ inherit the conformal invariance and 
restriction properties of the Brownian loop measure.

Now suppose that $\eta:[0,\infty) \rightarrow
\C$ is a simple curve
with $\eta(0,\infty) \subset \Half$
and $|\eta(t)|\rightarrow \infty$ as
$t \rightarrow \infty$.
Assume
that $\eta$ is parametrized by capacity, i.e.,
that
$\hcap[\eta[0,t]] =2 t$.  Let $H_t = \Half \setminus \eta[0,t]$
and let $g_t$ be the unique conformal transformation
of $H_t$ onto $\Half$ such that $g_t(\eta(t)) = 0$
and $g_t(z)\sim z $ 
 as $z \rightarrow \infty$. 
We let $f_t = g_t^{-1}$ which maps $\Half$ conformally
onto $H_t$ with $f_t(0) = \eta(t)$.

Given a realization ${\cal U}$ of the bubble
soup, consider the set of loops
\[   {\cal U}_{\eta,t} = \{f_{s_j} \circ \gamma_j:
   (\gamma_j,s_j) \in {\cal U}, s_j \leq t \} . \]
We consider this as realization of {\em unrooted}
loops by forgetting the loop.

\begin{theorem}  
\label{maintheorem}
For every $t < \infty$, if ${\cal U}$ is a bubble soup with intensity
$\lambda > 0$, then  ${\cal U}_{\eta,t}$, considered
as a collection of unrooted loops, is a realization of
${\cal L}^\perp(H_t)$
with intensity $\lambda$.
\end{theorem}

It is useful to consider this theorem in the
other direction, i.e. to see that it is equivalent to 
Theorem \ref{ls=bs}.  Let ${\cal L}$ be
a realization of the loop soup in $\Half$
with intensity $\lambda$.
We write elements of ${\cal L}$ as $[\gamma]$ since they are  
equivalence classes of loops. We write $V_\gamma$ for the
hull generated by $[\gamma]$, i.e. $V_\gamma$ is the complement of the unbounded
component of $\C \setminus \gamma[0,t_{\gamma}]$ (this does 
not depend on the choice of representative of $[\gamma]$). 
Let $\eta$ be as before
and let us 
write ${\cal L}^\perp = \{[\gamma_1],[\gamma_2],\ldots\}$ 
for ${\cal L}^\perp(\Half \setminus \eta[0,\infty))$, i.e.,
 for the set of loops in ${\cal L}$ that intersect
$\eta[0,\infty)$.   
For every $[\gamma_j] \in {\cal L} $,
let $r_j$ denote the smallest $r$ such that
$\eta(r) \in \gamma[0,t_\gamma]$. Note that this does
not depend on which representative $\gamma_j$ of
$[\gamma_j]$ that we choose. 

Let us now briefly justify the fact that 
   with probability
one, for each $j$ there is a unique
representative of $[\gamma_j]$, which we write
as just $\gamma_j$,  such that
$\gamma_j(0) = \eta(r_j)$ and $\gamma_j(0,t_{\gamma_j})
\subset H_{r_j}$. 
It follows for instance readily from the fact that if $B$ is a 
Brownian bridge (from $z$ to $w$ in time $t$), conditioned to stay in $\H$, then
for each rational $0  < q_1 < q_2 < t$, if one defines the 
first time $s(q_1)$ at which $\eta$ hits $B[0,q_1]$, then
\[ \Prob [\{\, s(q_1) < \infty; \;  \eta(s(q_1)) \in  B[q_2,t] \, \}] = 0 , \]
since complex Brownian motion does not hit points.

From now on we consider $[\gamma_j]$ as a rooted loop
by choosing this representative $\gamma_j$.  
Note that this choice depends on $\eta$. 
The set of times $ {\cal T} = \{r_j: \gamma_j \in {\cal L}^\perp \}$
is countable and   dense in $[0,\infty)$ since with
probability one
for each  rational $t$ there exists loops in ${\cal L}$ of
arbitrarily small diameter surrounding $\eta(t)$.
Also,  $r_j \neq r_k$ if $j \neq k$.  We let ${\cal L}^\perp_t$
denote the set of $\gamma_j \in {\cal L}$ with
$r_j \leq t$, i.e., the set of loops that intersect $\eta[0,t]$.  Recall
 that if $t < t_1$, then ${\cal L}^\perp_t$ and ${\cal L}^\perp_{t_1} \setminus {\cal L}^\perp_t$ are independent.
%
%
\medbreak

\Pf
  For $r > 0$, let ${\cal L}_t^\perp(r)$ denote the set of $\gamma_j
\in {\cal L}_t^\perp$ such that $\rad[g_{r_j} \circ \gamma_j] := \sup\{g_{r_j}
 \circ \gamma_j(s): 0 \leq s \leq t_{g_{r_j} \circ \gamma_j}\} \geq r$.  
Note that with probability one ${\cal L}_t(r)$ is finite for each $t < \infty, r > 0$. 
 It suffices to show that for every $r > 0$ the set of loops 
\[  \{g_{r_j} \circ \gamma_j: \gamma_j \in {\cal L}_t^\perp(r)\} \] is a
 Poissonian realization of the measure
$ \lambda t \, \bubble_\Half(0;r)$  We only need to do this for
 the case $r=1$; the other cases are essentially
the same.  Let ${\cal A}_t = \{g_{r_j} 
\circ \gamma_j: \gamma_j \in {\cal L}_t(1)^\perp\}$. 
We have already noted that for $\epsilon > 0$,
${\cal A}_{t + \epsilon}\setminus {\cal A}_t$ is independent of
${\cal A}_t$.   
If $t > 0$, the curve
 $\eta^t(s) = g_t[\eta(t+s)], 0 \leq s < \infty$, 
  is
also a simple
curve parametrized by capacity.   Conformal
invariance of $\loopmeasure$ tells us that the
distribution of $g_t \circ[{\cal A}_{t + \epsilon}
 \setminus {\cal A}_t]$,
derived from the curve $\eta$, is the same as the
distribution of ${\cal A}_\epsilon$ derived from the
curve $\eta^t$.  Hence it suffices to prove the
two conditions above for $t=0$.
But this is the estimate that was done
in \S\ref{estimatesec} so we have the result. 
\endpf

This implies (with (\ref {schw})) in particular immediately the following fact:
Suppose that $D \subset \H$ is simply connected, and that the curve 
$\eta (0,T] \subset D$ is parametrized as before. Define $D_t= g_t (D)$
(where $g_t$ is the conformal map from $\H \setminus \eta [0,t]$ onto $\H$ 
with $g_t (z) \sim z$ at infinity, and $g_t (\eta_t) =  O$).
As in (\ref {schw}), define also a conformal map $\phi_t$ from $D_t$ onto 
$\H$ that fixes the origin. Then, 
\begin {equation}
\label {schw2}
\Prob [ \forall \gamma \in {\cal L}\ : \ \gamma \cap \eta [0,T] =
\emptyset \hbox { or } \gamma \subset D ]
= \exp ( \lambda \int_0^T \frac {S_{\phi_t} (O)} 6 dt ).
\end {equation}

\subsection{Parametrization}

Suppose $\eta:[0,\infty) \rightarrow \C$ is a curve as before,  
 and let ${\cal L}^\perp_t$
denote a realization of the loop soup in $\Half$
restricted to curves
that intersect $\eta[0,t]$.   
We are going to show that if the path $\eta[0,\infty)$ has
dimension strictly less than two, then the sum of all the
time-lengths of the loops 
in ${\cal L}^\perp_t$ is almost surely finite. This will imply that one can 
construct a continuous path by attaching these loops
``chronologically'' to $\eta$.

\begin{lemma}  \label{parameterlemma}
Suppose that for some $\epsilon > 0$ and $T>0$,
\begin{equation}  \label{parameter2}
 \lim_{\delta \rightarrow 0+}
           \delta^{-\epsilon} \; 
  \area(\{z: \dist[z,\eta[0,T]] \leq
   \delta\})  = 0 . 
\end{equation}
Then with probability one, 
\[   \sum_{\gamma \in {\cal L}^\perp_{T}}
     t_\gamma < \infty . \]
\end{lemma}

\Pf Fix $T, \epsilon$, and let $r = \rad(\eta[0,T]) <
\infty$. Constants in this
proof may depend on $T ,r, \eps$. 
It suffices to prove two facts:  
\[ \#\{\gamma \in {\cal L}^\perp_{T}:
t_\gamma > 1\}  < \infty \hbox { a.s.,  and }\;\;\;
   \E[ \sum_{\gamma \in {\cal L}^\perp_{T}}
     t_\gamma  \; 1_{t_\gamma \leq 1}] < \infty . \]
Note that the first one is equivalent to 
$$ \loopmeasure_\Half [\{\gamma \in {\cal L}_T^\perp \ : \ t_\gamma > 1 \} ] 
< \infty .$$
But on the one hand 
\begin {eqnarray*}
\loopmeasure [ \{ \gamma \ : \ \gamma \subset 2r \Disk \ : \ t_\gamma > 1 \} ]
 &=&  \int_{2r \Disk} \int_1^\infty dA (z) \frac {dt}{2 \pi t^2}
\mu^\# (z,z; t) [\{ \gamma \ : \  \gamma \subset 2r \Disk \} ]\\ 
&\le& \frac {A (2r \Disk)}{2 \pi}
< \infty .
\end {eqnarray*}
On the other hand, 
\begin {eqnarray*}
\lefteqn {
\loopmeasure_\Half [ \{ \gamma \ : \ t_\gamma > 1 ,\  \gamma \not\subset 2r \Disk 
,\ \gamma \cap r\Disk \not= \emptyset \} ] 
} \\
&\le& 
\frac {1}{2 \pi} \int_0^r \int_0^\pi \bubble_{u\Disk_+} (u\exp (i \theta))
[ \{\gamma \ : \  \gamma \not\subset 2r \Disk \} ] du\  d\theta
.
\end {eqnarray*}
It is easy (using conformal invariance) to see that 
  $\bubble_{u\Disk_+} (u\exp (i \theta))
[ \{ \gamma \not\subset 2r \Disk \} ]$ is bounded independently from $u \le r$ and 
$\theta \in [0, \pi]$.
Hence, the last displayed expression is finite, which completes the proof of the fact
that the number of loops in ${\cal L}$
 of length greater than one and that do intersect $\eta [0,T]$
is almost surely finite.

Note that 
\[     \E[ \sum_{\gamma \in {\cal L}^\perp_{T}}
     t_\gamma  \; 1_{t_\gamma \leq 1}] 
  = \int_\Half \int_0^1 |\tilde \mu_\Half(z,z;t)|  \, dt \; dA(z), \]
where $\tilde \mu_\Half(z,z;t)$  denotes  $\mu_\Half(z,z;t)$
restricted to loops that intersect $\eta[0,T]$  (the $t_\gamma^{-1}$ in 
the definition of the loop measure cancels with the $t_\gamma$ in
the expression on the left hand side).   Let
\[    F(z) = \int_0^1   |\tilde \mu_\Half(z,z;t)|
  \; dt , \]
 and let $d_z = \dist(z,\eta[0,T])$.  It is standard to see that there exist
constants $c,a$ such that 
\[     |\tilde \mu_\Half(z,z;t)| \leq c \,t^{-1} \,
              e^{-a \, d_z^2/t} . \]
Hence, we get $F(z) \leq c \; \log (1/d_z)$ and
\begin{equation}  \label{parameter3}
\area\{z: F(z) \geq  s\}  \leq 
   \area\{z: \dist(z,\eta[0,T]) \leq  e^{-s/c} \}
   \leq e^{-s \epsilon/ c} . 
\end{equation}
Also we get $F(z) \leq c e^{-a |z|^2}$ for
$|z| \geq 3r$,  
and hence  we can see that $\int F(z) \; dA(z)
< \infty$. \endpf

\medskip

\noindent {\bf Remark.}  From (\ref{parameter3}) we can see that
(\ref{parameter2}) can be weakened to
\[  {\rm area}\{z: \dist(z,\eta[0,T]) \leq e^{-s} \} \leq g(s) , \]
where $\int_1^\infty g(s) \, ds < \infty$.  
However, if $\eta$ is space-filling, then the result does not hold as the 
following shows:

\begin {proposition}
If  $D$ is any nonempty open domain, then
$\sum_{\gamma \in {\cal L}(D)} t_\gamma = \infty$
almost surely.
\end {proposition}
 
\Pf
Note first that
\[   \E[\sum_{\gamma \in {\cal L}(D)} t_\gamma] = \infty . \]
This can be seen easily from the scaling rule
\[   \E[\sum_{\gamma \in {\cal L}(rD )} t_\gamma] 
    = r^2 \, \E[\sum_{\gamma \in {\cal L}(D)} t_\gamma]   . \]
For example, 
if  $D$ is a square, we can divide $D$ into $4$ squares of half
the side length, $D_1,\ldots,D_4$.   The scaling rule tells us
that
\begin{equation}  \label{parameter4}
  \E[\sum_{\gamma \in {\cal L}(D)} t_\gamma] = 
        \sum_{j=1}^4  \; \E[\sum_{\gamma \in {\cal L}(D_j)} t_\gamma] . 
\end{equation}
But
\[ \sum_{\gamma \in {\cal L}(D)} t_\gamma  =
     \sum_{j=1}^4  \sum_{\gamma \in {\cal L}(D_j)} t_\gamma  +
   \sum_{{\cal L}
  (D) \setminus [{\cal L}(D_1) \cup \cdots \cup {\cal L}(D_4)]}t_\gamma. \]
Since the last term has strictly positive expectation, the expectations
in (\ref{parameter4}) must be infinite. 

Furthermore (by dividing the square into 
$2^m$ smaller squares),
$\sum_{\gamma \in {\cal L}(D)} t_\gamma$ is larger than the mean of the values 
of $2^m$ independent copies of itself (i.e. the same random variable with infinite 
expectation) for any $m$. The result follows.
\endpf

\medskip

 With Lemma \ref{parameterlemma}
we can give a Brownian parametrization
to the curve ``$\eta$ with the loops added.''    Let
${\cal L}$ be a realization of the Brownian loop soup,
and let 
$\{[\gamma_1],[\gamma_2],\ldots\}$ be the  (unrooted) loops
that intersect $\eta[0,\infty)$.   As before, choose $r_j$ and
representative $\gamma_j$  so that $\gamma_j(0) =\eta(r_j)$
and $\gamma_j[0,t_{\gamma_j}] \cap \eta[0,r_j) = \emptyset$.
Define
\[       S(r-) = \sum_{r_j < r} t_{\gamma_j},
\;\;\;\;\;  S(r+) = \sum_{r_j \leq  r} t_{\gamma_j}. \]
Then $S(r)$ is an increasing function with jumps at
$r_j$ of size $t_{\gamma_j}$.  Define the process
$Y_s$ by
\[             Y_{S(r-)} = \eta(r) , \]
and if $S(r-) < S(r+)$,
\[           Y_{S(r-) + s} = \gamma_j(s),
 \;\;\;\;  0 \leq s \leq t_\gamma . \]
The density of the loop-soup implies readily that $t \mapsto Y_t$
is continuous (provided $\eta$ is a simple curve for instance).

The results of \cite {LSWlesl} strongly suggest that the following conjecture holds.

\begin {conjecture}
If the curve $\eta$ is chordal $SLE_2$, and 
$\lambda = 1$, then the law of $Y$ is $\mu_\Half^\normed (0, \infty)$.
\end {conjecture} 

There seem to be different possible ways to prove this. One can
use the convergence of loop-erased random walk to the $SLE_2$ curve \cite {LSWlesl}. The main missing step is the convergence of discrete bubbles towards the Brownian bubbles.
 
\medbreak

\noindent 
{\bf Acknowledgements.}
We would like to thank Oded Schramm for many inspiring conversations.
Part of this work was carried at the Centre Emile Borel of the Institut Henri Poincar\'e.

Gregory Lawler

Department of Mathematics

310 Malott Hall

Cornell University

Ithaca, NY 14853-4201, USA 

lawler@math.cornell.edu

\medbreak

Wendelin Werner

Laboratoire de Math\'ematiques 

B\^at. 425

Universit\'e Paris-Sud

91405 Orsay cedex, France

wendelin.werner@math.u-psud.fr

\end{document}